\documentclass[11pt]{amsart}
\usepackage{amscd}
\usepackage{amssymb}
\usepackage{hyperref}

\newtheorem{thm}{Theorem}[section]

\theoremstyle{definition}

\renewcommand{\thecase}{}

\newtheorem{conj}[thm]{Conjecture}

\renewcommand{\thestep}{}

\theoremstyle{remark}

\makeatletter
\def\alphenumi{
  \def\theenumi{\alph{enumi}}
  \def\p@enumi{\theenumi}
  \def\labelenumi{(\@alph\c@enumi)}}
\makeatother

\makeatletter
\def\thecase{\@arabic\c@case}
\makeatother

\makeatletter
\def\thestep{\@arabic\c@step}
\makeatother

\newcommand\embed{\hookrightarrow}

\newcommand\AAA{\mathbb{A}}

\newcommand\CC{\mathbb{C}}

\newcommand\NN{\mathbb{N}}

\newcommand\QQ{\mathbb{Q}}
\newcommand\RR{\mathbb{R}}

\newcommand\ZZ{\mathbb{Z}}

\newcommand\bgamma{{\boldsymbol{\gamma}}}
\newcommand\bchi{{\boldsymbol{\chi}}}

\newcommand\bD{{\mathbf{D}}}

\newcommand\bL{{\mathbf{L}}}

\newcommand\bS{{\mathbf{S}}}

\newcommand\bW{{\mathbf{W}}}
\newcommand\bx{{\mathbf{x}}}

\newcommand{\rd}{\partial}

\newcommand\thalf{{\textstyle{\frac{1}{2}}}}

\newcommand\half{{{\frac{1}{2}}}}

\newcommand\quarter{{{\frac{1}{4}}}}

\newcommand\fg{{\mathfrak{g}}}

\newcommand\fs{{\mathfrak{s}}}

\newcommand\ft{{\mathfrak{t}}}

\newcommand\eps{\varepsilon}

\newcommand\La{\Lambda}

\newcommand\si{\sigma}
\newcommand\Si{\Sigma}

\newcommand\su{{\mathfrak{s}\mathfrak{u}}}

\newcommand\PU{\operatorname{PU}}

\newcommand\SO{\operatorname{SO}}

\newcommand\less{\setminus}

\newcommand{\8}{\infty}

\newcommand\ad{{\operatorname{ad}}}

\newcommand\CCl{\operatorname{{\mathbb{C}\ell}}}

\newcommand\End{\operatorname{End}}

\newcommand\Fr{\operatorname{Fr}}

\newcommand\Gl{\operatorname{Gl}}

\newcommand\Hom{\operatorname{Hom}}

\newcommand\Map{\operatorname{Map}}

\newcommand\Sym{\operatorname{Sym}}

\newcommand\Tr{\operatorname{Tr}}

\newcommand\vol{\operatorname{vol}}

\newcommand\spinc{\text{$\text{spin}^c$ }}
\newcommand\spinu{\text{$\text{spin}^u$ }}
\newcommand\Spinc{\text{$\text{Spin}^c$}}

\newcommand\sG{{\mathcal{G}}}

\newcommand\sM{{\mathcal{M}}}

\newcommand\tN{{\tilde N}}

\begin{document}
\title[On Donaldson and Seiberg-Witten invariants]
{On Donaldson and Seiberg-Witten invariants}
\author[Paul M. N. Feehan]{Paul M. N. Feehan}
\address{Department of Mathematics\\
Rutgers University\\
Piscataway, NJ 08854-8019}
\email{feehan@math.rutgers.edu}
\urladdr{http://www.math.rutgers.edu/$\sim$feehan}
\thanks{PMNF was supported in part by NSF grants DMS-9704174 and
DMS-0125170}

\author[Thomas G. Leness]{Thomas G. Leness}
\address{Department of Mathematics\\
Florida International University\\
Miami, FL 33199}
\email{lenesst@fiu.edu}
\urladdr{http://www.fiu.edu/$\sim$lenesst}
\thanks{TGL was supported in part by NSF grants DMS-0103677}
\dedicatory{}
\date{This version: August 7, 2002. arXiv:math.DG/0106221}
\begin{abstract}
  We sketch a proof of Witten's
  formula relating the Donaldson and Seiberg-Witten series modulo
  powers of degree $c+2$, with $c=-\frac{1}{4}(7\chi+11\sigma)$, for
  four-manifolds obeying some mild conditions, where $\chi$ and $\sigma$
  are their Euler characteristic and signature. We use the moduli space of
  $\SO(3)$ monopoles as a cobordism between a link of the Donaldson moduli
  space of anti-self-dual $\SO(3)$ connections and links of the moduli
  spaces of Seiberg-Witten monopoles. Gluing techniques allow us to compute
  contributions from Seiberg-Witten moduli spaces lying in the first (or
  `one-bubble') level of the Uhlenbeck compactification of the moduli
  space of $\SO(3)$ monopoles.
\end{abstract}
\maketitle


\newcommand\vir{{\mathrm{vir}}}

\section{Introduction}
This article consists of lightly edited notes for a lecture by the first
author at the International Georgia Topology Conference 2001. Although we
shall only briefly mention technical details and qualifications appropriate
for more complete accounts published elsewhere \cite{FL2a}, \cite{FL2b},
\cite{FLLevelOne}, we hope that these notes provide a convenient survey of our
recent work on the SO(3)-monopole program.

\subsection{Witten's conjecture}
Two kinds of invariants can be used to explore the classification problem for
compact, smooth 4-manifolds:

\begin{itemize}
\item
{\em Donaldson invariants\/}, defined using an SO(3) Yang-Mills gauge
theory (discovered in 1986).
\item
{\em Seiberg-Witten invariants\/}, defined using a U(1) monopole gauge theory
(1994).
\end{itemize}
We shall restrict our attention throughout to the case of closed, oriented
4-manifolds with $b_1=0$ and odd $b_2^+>1$.  The conjectured relationship
between these gauge theory invariants is described below:

\begin{conj}
\label{conj:Witten}
\cite{Witten}
A 4-manifold $X$ has KM-simple type if and only if it has SW-simple
type. If $X$ has simple type, the KM and SW basic classes coincide, and
the Donaldson and Seiberg-Witten series obey
\begin{equation}
\label{eq:WittenFormula}
\bD_X^w(h) = 2^{2-c(X)}e^{\frac{1}{2}Q_X(h,h)}\bS\bW_X^w(h),
\quad h\in H_2(X;\RR).
\end{equation}
\end{conj}

Here, the 4-manifold $X$ has intersection form
$$
Q_X:H_2(X;\ZZ)\times H_2(X;\ZZ)\to \ZZ,
$$
Euler characteristic $\chi$, signature $\sigma$ and
$$
c(X)=-\frac{1}{4}(7\chi+11\sigma).
$$
We shall recall the definitions of the Donaldson and Seiberg-Witten series
shortly.

\subsection{Remarks on the problem}
Before proceeding to discuss our work on Witten's conjecture, it is
interesting to compare the mathematicians' and physicists' approaches to
establishing \eqref{eq:WittenFormula}.

Witten employs a certain {\em $N=2$ supersymmetric quantum Yang-Mills
theory\/}. He uses rescaling, $g_t=t^2g$, of the Riemannian metric $g$ on $X$
and metric independence of the correlation functions to relate the
Donaldson invariants ($t\to 0$) with the Seiberg-Witten invariants ($t\to
\8$). 

The mathematical approach to a proof of Witten's formula proposed by
Pidstrigatch and Tyurin \cite{PTLocal} instead employs
an SO(3) monopole gauge theory which generalizes both the instanton and
U(1) monopole gauge theories. All three gauge theories are classical field
theories and their solutions are invariant under metric rescaling, whereas
Witten's quantum field theory is sensitive to metric rescaling.

Apparently, the SO(3) monopole gauge theory provides a purely classical
field theory alternative to Witten's quantum field theory method. The
problem of determining the relationship between these two approaches is
surely an important one worth exploring further.

\section{SO(3) monopoles}
\subsection{Clifford modules and spin structures}
Given a Riemannian metric $g$ on $X$, let $V\to X$ be a Hermitian bundle with
a linear Clifford map $\rho:T^*X\to\End_\CC(V)$,
$$
\rho(\alpha)^\dagger = -\rho(\alpha),
\quad\text{and}\quad
\rho(\alpha)^2 = -g(\alpha,\alpha),\quad \alpha\in\Omega^1(X,\RR).
$$
Then $(\rho,V)$ defines a Clifford or $\CCl(T^*X)$ module structure on
$V$.  If $W\to X$ is a complex-rank four Hermitian bundle, then 
$\fs=(\rho,W)$ is a {\em \spinc structure\/}, familiar from Seiberg-Witten
theory \cite{SalamonSWBook}. If $V\to X$ is a complex-rank eight Hermitian
bundle, then we call $\ft=(\rho,V)$ a {\em \spinu structure\/}.

\subsection{From \spinu\  structures to SO(3) bundles}
Given $\ft=(\rho,V)$ on $X$, one obtains
\begin{itemize}
\item
An {\em $\SO(3)$ subbundle\/},
$$
\fg_{\ft} \subset \su(V),
$$
characterized as the span of the sections $\xi$ of $\su(V)$ such that
$[\xi,\rho(\omega)]=0$, for all $\omega\in\Omega^\bullet(X,\RR)$.
\item
A complex line bundle, 
$$
{\det}^{\frac{1}{2}}(V^+),
$$
where $V = V^+\oplus V^-$ and $V^\pm$ are the $\mp 1$ eigenspaces of
$\rho(\vol)$ on $V$.
\item
Splittings, if $\Lambda^2(T^*X) = \Lambda^+\oplus\Lambda^-$ and
$\Lambda^\pm$ are the $\pm 1$ eigenspaces of $*_g$ on $\Lambda^2(T^*X)$,
and $\rho:\Lambda^\pm\cong\su(W^\pm)$ are the usual isomorphisms of SO(3)
bundles, 
$$
\su(V^\pm)
\cong
\rho(\Lambda^\pm)\oplus\rho(\Lambda^\pm)\otimes\fg_{\ft}\oplus\fg_{\ft}.
$$
\end{itemize}

Moreover, for any choice of \spinc structure $\fs=(\rho,W)$, one further
obtains 
\begin{itemize}
\item
A complex-rank two Hermitian bundle, 
$$
E = \Hom_{\CCl(T^*X)}(W,V).
$$
\item
A Clifford module isomorphism, 
$$
V \cong W\otimes E.
$$
\item 
An isomorphism of SO(3) bundles
$$
\su(E)\cong \fg_{\ft}.
$$
\item
An isomorphism of complex line bundles,
$$
{\det}^{\frac{1}{2}}(V^+) \cong \det(W^+)\otimes\det(E).
$$
\end{itemize}

\subsection{SO(3)-monopole equations} 
We call a pair $(A,\Phi)$ an {\em $\SO(3)$ monopole\/} if
\begin{align*}
\ad^{-1}(F_{\hat A}^+) - \rho^{-1}(\Phi\otimes\Phi^*)_{00}
&=0,
\\
D_A\Phi 
&=0,
\end{align*}
where $A$ is a spin connection on $V$, inducing a fixed connection
$A|_{\det(V^+)}=2A_\Lambda$ on $\det(V^+)$ and $\Phi$ is a section of $V^+$;
$\hat A$ is the induced connection on the $\SO(3)$ bundle
$\fg_{\ft}\subset\su(V)$; 
$F_{\hat A}^+$ is the self-dual component of the curvature of $\hat A$; the
term $(\Phi\otimes\Phi^*)_{00}$ is the component of $\Phi\otimes\Phi^*$
lying in $\rho(\Lambda^+)\otimes\fg_{\ft}$;
$D_A:C^\8(V^+)\to C^\8(V^-)$ is the Dirac operator.
We let $\sM_{\ft}$ be the space of SO(3) monopoles for $\ft=(\rho,V)$,
modulo gauge transformations.

\subsection{Singularities in SO(3)-monopole space}
We now classify the fixed points of the circle action on $\sM_{\ft}$ given
by complex multiplication on the spinor components.

An SO(3) monopole 
$(A,\Phi)$ is a {\em Yang-Mills\/} or {\em instanton\/} solution if
$$
\Phi\equiv 0 \quad \text{and}\quad F_{\hat A}^+ = 0.
$$
Hence, there is a moduli subspace of SO(3) instantons,
$$
M_\kappa^w \embed \sM_{\ft},
$$
where $\kappa=-\quarter p_1(\fg_{\ft})$ and $w\in H^2(X;\ZZ)$ lifts
$w_2(\fg_{\ft})\in H^2(X;\ZZ/2\ZZ)$. 

An SO(3) monopole 
$(A,\Phi)$ is a {\em Seiberg-Witten\/} or {\em reducible\/} solution if
$$
A = B\oplus B\otimes A_L 
\quad\text{on}\quad
V = W\oplus W\otimes L,
$$
for some Hermitian line bundle $L$, a unitary connection $A_L$ on $L$,
and $\Phi = \Psi\oplus 0$ with $\Psi$ a section of $W^+$ obeying
\begin{align*}
\Tr(F_{B}^+) - \rho^{-1}(\Psi\otimes\Psi^*)_0 - F_{A_\Lambda}^+ &= 0, 
\\
D_B\Psi &= 0.
\end{align*}
Hence, there are moduli subspaces of Seiberg-Witten monopoles for
$\fs=(\rho,W)$, 
$$
M_{\fs} \embed \sM_{\ft},
$$
whenever $V=W\oplus W\otimes L$.

\section{Invariants of smooth 4-manifolds}
We sketch definitions of the Donaldson series
\cite{KMStructure} and Seiberg-Witten series \cite{Witten}.

\subsection{Donaldson invariants}
Set $\AAA(X) = \Sym(H_0(X;\RR)\oplus H_2(X;\RR))$, so $z\in\AAA(X)$ is a
linear combination of monomials
$$
x^m\beta_1\beta_2\cdots\beta_{\delta-2m},
$$
with $x\in H_0(X;\ZZ)$ being the positive generator and $\beta_i \in
H_2(X;\RR)$.  Cohomology classes on $M_\kappa^w$ can be defined via a map
\cite{DonPoly}, \cite{DK}, 
$$
\mu_p:H_i(X;\RR) \to H^{4-i}(M_\kappa^w;\RR).
$$
The Donaldson invariant is then a linear function
$$
D_X^w:\AAA(X)\to\RR,
$$
where, for a monomial $z$ with $\deg(z)=2\delta$,
$$
D_X^w(z) = \langle \mu_p(z),[M_\kappa^w]\rangle
$$
with $\mu_p(z) =
\mu_p(x)^m\smile\mu_p(\beta_1)\smile\cdots\smile\mu_p(\beta_{\delta-2m})$.

\subsection{Kronheimer-Mrowka structure theorem}
One says that a 4-manifold $X$ has {\em KM-simple type\/} if for some $w$
and all $z\in\AAA(X)$, 
$$
D_X^w(x^2z) = 4D_X^w(z).
$$
One defines the {\em Donaldson series\/} by setting
\begin{equation}
\label{eq:DSeries}
\bD_X^w(h) = D_X^w((1+\thalf x)e^h), \quad h\in H_2(X;\RR).
\end{equation}
We recall the celebrated

\begin{thm}
\cite{KMStructure}
If $X$ has KM-simple type, then there exist $a_r\in\QQ$ and $K_r\in
H^2(X;\ZZ)$, the KM-basic classes, such that
\begin{equation}
\label{eq:KMStructure}
\bD^w_X(h) 
= 
e^{\frac{1}{2}Q_X(h,h)}
\sum_{r=1}^s(-1)^{\frac{1}{2}(w^2+w\cdot K_r)}a_r e^{\langle K_r,h\rangle},
\quad h\in H_2(X;\RR).
\end{equation}
\end{thm}

\noindent See also \cite{FSStructure}, for a similar result and independent
proof by different methods.

\subsection{Seiberg-Witten invariants}
The {\em Seiberg-Witten invariants\/} comprise a function, 
$$
SW_X:\Spinc(X) \to \ZZ,
$$
where
$$
SW_X(\fs) = \langle\mu_{\fs}(x)^{\max},[M_{\fs}]\rangle,
$$
and $\mu_{\fs}(x)\in H^2(M_{\fs};\ZZ)$ is a cohomology class associated to
a circle action.
One says that a 4-manifold $X$ has {\em SW-simple type\/} if for all $\fs$
for which $M_{\fs}$ has positive dimension one has that
$$
SW_X(\fs) = 0,
$$
and calls $c_1(\fs) = c_1(W^+)$ an
{\em SW-basic class} if $SW_X(\fs)\neq 0$.
We define the {\em Seiberg-Witten series\/} by
\begin{equation}
\label{eq:SWSeries}
\bS\bW_X^{w}(h) =\sum_{\fs \in
\Spinc(X)}(-1)^{\half(w^{2}+c_{1}(\fs)\cdot w)}
SW_X(\fs)e^{\langle c_{1}(\fs),h\rangle},
\end{equation}
for all $h\in H_2(X;\RR)$.
Witten's prediction then takes the form stated in Conjecture
\ref{conj:Witten}.

\section{SO(3) monopole cobordism}
\subsection{Bubbling and Uhlenbeck compactness}
In order to apply the moduli space $\sM_{\ft}$ of SO(3) monopoles as a
cobordism, we must use a compactification.
If $\{[A_\alpha,\Phi_\alpha]\}_{\alpha\in\NN}\subset \sM_{\ft}$, then
the sequence converges to an {\em ideal $\SO(3)$ monopole\/}
$([A_\8,\Phi_\8],\bx)$ in $M_{\ft_{\ell}}\times\Sym^\ell(X)$ if
\begin{itemize}
\item 
$(A_\alpha,\Phi_\alpha) \to (A_\8,\Phi_\8)$ in $C^\8$ on $X\less\{\bx\}$,
modulo gauge transformations,
\item
$|F_{A_\alpha}|^2 \to |F_{A_\8}|^2 + 8\pi^2\sum_{x\in\bx}\delta(x)$, as
measures, where $\delta(x)$ denotes the Dirac measure centered at $x$. 
\end{itemize}
We let $\bar\sM_{\ft}$ be the closure of $\sM_{\ft}$ with respect to the Uhlenbeck
topology, implicit above, in the space of ideal SO(3) monopoles,
$$
\bigsqcup_{\ell = 0}^N \left(M_{\ft_{\ell}}\times\Sym^\ell(X)\right),
$$
where $\ft_{\ell}=(\rho,V_{\ell})$ is a \spinu structure with
characteristic classes
$$
p_1(\fg_{\ft_\ell})=p_1(\fg_{\ft})+4\ell,
\quad
w_2(\fg_{\ft_{\ell}}) = w_2(\fg_{\ft}),
\text{ and }
c_1(\ft_{\ell}) = c_1(\ft), 
$$
and $c_1(\ft) := c_1({\det}^{\frac{1}{2}}(V^+))$.  The space
$\bar\sM_{\ft}$ is smoothly stratified, with {\em top\/} or {\em zeroth
level\/} $\sM_{\ft}$, and {\em lower levels\/}
$\sM_{\ft_{\ell}}\times\Sym^\ell(X)$, $\ell \geq 1$.

\subsection{Stratification of the space of SO(3) monopoles}
For the top level, one has a stratification
\begin{equation}
\label{eq:TopStratification}
\sM_{\ft}
=
M_\kappa^w
\sqcup
\sM_{\ft}^{*,0}
\sqcup
\bigcup_{\ft=\fs\oplus\fs'}
M_{\fs}.
\end{equation}
The complement $\sM_{\ft}^{*,0}$ in $\sM_{\ft}$ of the Yang-Mills and
Seiberg-Witten solutions is a smooth manifold, cut out transversely by the
SO(3) monopole equations \cite{FeehanGenericMetric}. 

A stratification of the form \eqref{eq:TopStratification} arises in each
level, $\sM_{\ft_\ell}\times\Sym^\ell(X)$, of the compactification
$\bar\sM_{\ft}$. Though dimension-counting arguments rule out contributions
from the instanton moduli subspace $M_{\kappa-\ell}^w$ of $\sM_{\ft_\ell}$
to pairings with the cohomology classes appearing in equation
\eqref{eq:RawIdentity}, Seiberg-Witten moduli subspaces of $\sM_{\ft_\ell}$
can contribute to Donaldson invariants computed using
$M_\kappa^w\subset\sM_{\ft}$. To apply the cobordism, we
\begin{itemize}
\item
define a link $\bar\bL_{\ft,\kappa}^w\subset \bar\sM_{\ft}/S^1$ of the
instanton moduli subspace, 
$$
\bar M_\kappa^w \subset \bar\sM_{\ft},
$$
by restricting to spinors with $L^2$ norm equal to a small positive
constant, and 
\item
use gluing theory \cite{FL3}, \cite{FL4} to construct links 
$\bar\bL_{\ft,\fs}\subset \bar\sM_{\ft}/S^1$ of ideal
Seiberg-Witten moduli subspaces, 
$$
M_{\fs}\times\Sym^\ell(X) \subset \bar\sM_{\ft}.
$$
The links $\bar{\bL}_{\ft,\fs}$ are considerably more difficult to
construct than $\bar\bL_{\ft,\kappa}^w$, especially when $\ell$ is large.
\end{itemize}

\subsection{SO(3)-monopole cobordism formula}
The cobordism $\bar\sM_{\ft}/S^1$ now yields the raw identity,
\begin{equation}
\label{eq:RawIdentity}
\left\langle\mu_{p}(z)\smile \mu_{c}^{\delta_c-1},
[\bar\bL_{\ft,\kappa}^w]\right\rangle
=
-\sum_{\fs\in\Spinc(X)}
\left\langle\mu_{p}(z)\smile \mu_{c}^{\delta_c-1},
[\bar\bL_{\ft,\fs}]\right\rangle,
\end{equation}
where $\mu_c \in H^2(\sM_{\ft}^{*,0};\ZZ)$ is a class associated to a
circle action on $\sM_{\ft}^{*,0}$. One finds that
$$
\langle\mu_{p}(z)\smile \mu_{c}^{\delta_c-1},[\bar\bL_{\ft,\kappa}^w]\rangle
$$ 
is a multiple of the Donaldson invariant, $D_X^w(z)$.  The sum in
\eqref{eq:RawIdentity} is over all $\fs\in\Spinc(X)$, with $\bar\bL_{\ft,\fs}$
empty unless $M_{\fs}\times\Sym^\ell(X)\subset\bar\sM_{\ft}$, for some
$\ell(\ft,\fs)\geq 0$. 

The difficult aspect of using \eqref{eq:RawIdentity} to derive Witten's formula
\eqref{eq:WittenFormula} is to show that
\begin{equation}
\label{eq:SWLinkPairing}
\langle\mu_{p}(z)\smile \mu_{c}^{\delta_c-1}, [\bar\bL_{\ft,\fs}]\rangle
\end{equation}
is the correct multiple of the Seiberg-Witten invariant $SW_X(\fs)$; the
degree of difficulty grows rapidly with $\ell\geq 0$.
The assertion that the pairing \eqref{eq:SWLinkPairing} is a multiple
of $SW_X(\fs)$ is refered to as the  multiplicity conjecture.
As this conjecture follows from the work in \cite{FL5}, we shall
assume it for the rest of this note.

\section{Application of the cobordism}
We may consider following situations, arranged in increasing order of
complexity:
\begin{itemize}
\item 
There are no Seiberg-Witten moduli spaces with non-zero invariants in
$\bar\sM_{\ft}$, so the intersection 
$(M_{\fs}\times\Sym^\ell(X))\cap\bar\sM_{\ft}$
is empty for all $\ell\geq 0$ and $SW_X(\fs)\neq 0$. The Donaldson
invariants defined by $\bar M_\kappa^w\subset\bar\sM_{\ft}$ are then zero
and, eventually, this observation leads to a vanishing result:
$$
\bD_X^w(h) \equiv 0 \equiv \bS\bW_X^w(h) \pmod{h^{c(X)-2}},
\quad 
h\in H_2(X;\RR).
$$
\item
Calculation of contributions from $M_{\fs}\subset\sM_{\ft}$ leads to
Witten's formula, mod $h^{c(X)}$.
\item
Calculation of contributions from 
$M_{\fs}\times\Sym^\ell(X)\subset\bar\sM_{\ft}$ for $\ell =0,1$ leads to
Witten's formula, mod $h^{c(X)+2}$.
\item
Calculation of contributions from
$M_{\fs}\times\Sym^\ell(X)\subset\bar\sM_{\ft}$ for $\ell =0,1,2$ leads to
Witten's formula, mod $h^{c(X)+4}$.
\item
Calculation of contributions from
$M_{\fs}\times\Sym^\ell(X)\subset\bar\sM_{\ft}$ for $\ell \geq 3$ should
lead to a verification of Witten's formula \eqref{eq:WittenFormula}.
\end{itemize}
We have considered the cases $\ell=0$ and $\ell=1$ in detail \cite{FL2a},
\cite{FL2b}, \cite{FLLevelOne} and we would expect the case $\ell=2$ to follow
in a similar manner, by exploiting work of Leness \cite{LenessWC} on
the wall-crossing formula for Donaldson invariants.

At present, we can compute the general shape of the contributions for
$\ell\geq 3$ (see \cite{FL5}); complete, direct computations of those
contributions appear to be difficult, though we expect indirect methods
will yield the desired result \cite{FLWConjecture}.

\subsection{Level-zero Seiberg-Witten contributions}
Let $B \subset H^2(X;\ZZ)$ be the set of Seiberg-Witten basic classes and
let $B^\perp\subset H^2(X;\ZZ)$ be the $Q_X$-orthogonal complement of $B$.
We call a 4-manifold $X$ {\em abundant\/} if $Q_X|_{B^\perp}$ has a
hyperbolic sublattice.  Every compact, complex algebraic, simply-connected
surface with $b_2^+ \geq 3$ is abundant \cite{FL2a}.

\begin{thm}
\label{thm:WCL0}
\cite{FL2b}
Assume $X$ is abundant, has $b_1=0$, odd $b_2^+\geq 3$, and SW-simple type.
Suppose $\Lambda \in B^\perp$ exists with $\Lambda^2=2-(\chi+\sigma)$. For
such $\Lambda$ and $w\in H^2(X;\ZZ)$ with $w-\Lambda\equiv w_2(X)\pmod{2}$
one has, for all $h\in H_2(X;\RR)$,
\begin{align}
\label{eq:MMP}
\bD^{w}_X(h) 
&\equiv 0 \equiv \bS\bW^{w}_X(h)
\pmod{h^{c(X)-2}},
\\
\label{eq:WCL0}
\bD^{w}_X(h) 
&\equiv 2^{2-c(X)}e^{\half Q_X(h,h)}\bS\bW^{w}_X(h)
\pmod{h^{c(X)}}.
\end{align}
\end{thm}

The vanishing assertion \eqref{eq:MMP} for the Seiberg-Witten series is a
statement that the Moore-Mari\~no-Peradze conjecture holds for (abundant)
4-manifolds of SW-simple type \cite{FKLM}, \cite{MMPdg},
\cite{MMPhep}.  

\subsection{Level-one Seiberg-Witten contributions}
With more sophisticated analytical tools, specifically gluing theory, we can
compute contributions from $M_{\fs}\times X\subset \bar\sM_{\ft}$, and these
computations lead to the

\begin{thm}
\label{thm:WCL1}
\cite{FLLevelOne}
Same hypotheses as Theorem \ref{thm:WCL0}, but now
suppose $\Lambda \in B^\perp$ exists with
$\Lambda^2=4-(\chi+\sigma)$. For such $\Lambda$ and $w\in H^2(X;\ZZ)$ with
$w-\Lambda\equiv w_2(X)\pmod{2}$ one has, for all $h\in H_2(X;\RR)$,
\begin{align}
\label{eq:MMPagain}
\bD^{w}_X(h) 
&\equiv 0 \equiv \bS\bW^{w}_X(h)
\pmod{h^{c(X)-2}},
\\
\label{eq:WCL1}
\bD^{w}_X(h) 
&\equiv 2^{2-c(X)}e^{\half Q_X(h,h)}\bS\bW^{w}_X(h)
\pmod{h^{c(X)+2}}.
\end{align}
\end{thm}

We expect an identity similar to \eqref{eq:WCL1}, but mod $h^{c(X)+4}$, by
computing contributions for $\ell=2$, when $\Lambda^2 = 6-(\chi+\sigma)$.
The restrictive hypotheses on existence of classes $\Lambda$ with
prescribed even squares can be dropped if one can consider contributions
for arbitrary $\ell\geq 0$.

\subsection{Seiberg-Witten contributions from arbitrary levels}
More generally, we establish the following in \cite{FL5}:

\begin{thm}
\label{thm:SharpPTConjecture}
Let $X$ be a closed, connected, oriented smooth four-manifold with
$b_1(X)=0$ and odd $b_2^+(X)>1$. Let $\La,w\in H^2(X;\ZZ)$ obey
$w-\La\equiv w_2(X)\pmod 2$. Let $\delta,m$ be non-negative integers for
which $m\leq [\delta/2]$, where $[\,\cdot\,]$ denotes the greatest integer
function, and $\delta\equiv -w^2-\frac{3}{4}(\chi+\sigma)\pmod{4}$, with
$\Lambda$ and $\delta$ obeying $\delta<i(\La)$, where
$i(\La)=\La^2-\frac{1}{4}(\chi+\si)$. Then for any $h\in H_2(X;\RR)$ and
generator $x\in H_0(X;\ZZ)$, we have the following expression for the
Donaldson invariant:
\begin{equation}
\label{eq:RoughWittenFormula}
\begin{aligned}
D^w_X(h^{\delta-2m}x^m)
&=
\sum_{\fs\in{\textrm{\em \Spinc}}(X)}
(-1)^{\half(w^2+w\cdot c_1(\fs))}SW_X(\fs)
\\
&\qquad\times\sum_{i=0}^{ \min(\ell,[\delta/2]-m)}
\left(
p_{\delta,\ell,m,i}(c_1(\fs)-\La,\La)Q_X^i
\right)(h),
\end{aligned}
\end{equation}
where $Q_X$ is the intersection form on $H_2(X;\RR)$,
$\ell=\frac{1}{4}(\delta+(c_1(\fs)-\La)^2+\frac{3}{4}(\chi+\si))$
and $p_{\delta,\ell,m,i}(\cdot,\cdot)$ is a 
homogeneous polynomial of degree $\delta-2m-2i$
with coefficients which are universal functions of
$$
\chi,\ \sigma,\ c_1(\fs)^2,\ \La^2,\ c_1(\fs)\cdot\La,\ \delta,\ m, \ \ell.
$$
\end{thm}

Although Theorem \ref{thm:SharpPTConjecture} does not immediately yield
Witten's formula \eqref{eq:WittenFormula}, it is still powerful enough to
prove that the Seiberg-Witten invariants determine the Donaldson
invariants. Furthermore, Witten's formula \eqref{eq:WittenFormula} itself
should follow from Theorem \ref{thm:SharpPTConjecture} by indirect
calculations of the remaining unknown coefficients \cite{FLWConjecture}.

\section{Outline of the proofs of Theorems \ref{thm:WCL0}, \ref{thm:WCL1},
  and \ref{thm:SharpPTConjecture}}
We shall first sketch how to compute the rough form of the pairings,
\begin{equation}
\label{eq:BasicPairings}
\langle 
\mu_{p}(z)\smile\mu_{c}^{\delta_c-1}, [\bar\bL_{\ft,\fs}]
\rangle,
\end{equation}
or, at least why these pairings have the form
$$
SW_X(\fs) \times \text{(Factors depending only on topology)}.
$$
We use our gluing theory \cite{FL3}, \cite{FL4} to construct a topological
model for a neighborhood in $\bar\sM_{\ft}$ and hence a link,
$\bar\bL_{\ft,\fs}\subset\bar\sM_{\ft}/S^1$, of the `stratum'
$$
M_{\fs}\times\Sym^\ell(X)\subset \bar\sM_{\ft}.
$$
Given this topological model for $\bar\bL_{\ft,\fs}$, we can then apply
intersection theory methods to partly compute the pairings
\eqref{eq:BasicPairings}. 

This suffices to prove the `rough version' 
\eqref{eq:RoughWittenFormula} of Witten's formula \eqref{eq:WittenFormula},
which is enough to show that the
Seiberg-Witten invariants determine the Donaldson invariants. We
shall illustrate the method below, often assuming $\ell=1$ for the sake of
simplicity \cite{FLLevelOne}.

The passage from this stage to Witten's formula \eqref{eq:WittenFormula}
requires us to compute the many universal, but unknown coefficients in the
rough version \eqref{eq:RoughWittenFormula} of Witten's formula.

\subsection{Neighborhood of a Seiberg-Witten stratum}
Our gluing theory \cite{FL3}, \cite{FL4} allows us to construct a model for
a neighborhood of the level $M_{\fs}\times\Sym^\ell(X)$ in a local, `virtual'
moduli space,
\begin{equation}
\label{eq:GluingModel}
\bar\sM_{\ft,\fs}^{\vir}
:=
\tN_{\ft_\ell,\fs}(\eps)\times_{\sG_{\fs}}\bar{\Gl}_{\ft_\ell}(\delta),
\end{equation}
where $\bar{\Gl}_{\ft_\ell}(\delta)$ is a space of instanton gluing data
and $\tN_{\ft_\ell,\fs}(\eps)$ is a radius-$\eps$ disk subbundle of the
vector bundle \eqref{eq:VirtualNormalBundle}.

A neighborhood of $M_{\fs}\times\Sym^\ell(X)$ in the true moduli space
$\bar\sM_{\ft}$ then takes the shape
$$
\bgamma\left(\bchi^{-1}(0)\cap\bar\sM_{\ft,\fs}^{\vir}\right)
=
\bar\sM_{\ft}\cap\bgamma(\bar\sM_{\ft,\fs}^{\vir}),
$$
where $\bgamma$ and $\bchi$ are described below.  The stratum $M_{\fs}$ in
$\sM_{\ft_\ell}$ has a `virtual' normal bundle,
\begin{equation}
\label{eq:VirtualNormalBundle}
N_{\ft,\fs}\to M_{\fs},
\end{equation}
and an obstruction bundle with section $\bchi$,
$$
\Xi_{\ft,\fs}\to \sM_{\ft_\ell},
$$
while the gluing map,
$$
\bgamma:\bar\sM_{\ft,\fs}^{\vir}
\to
\text{Configuration space of ideal pairs containing $\bar\sM_{\ft}$},
$$
gives a homeomorphism from $\bchi^{-1}(0)\cap\bar\sM_{\ft,\fs}^{\vir}$
onto a neighborhood of $M_{\fs}\times\Sym^\ell(X)$ in $\bar\sM_{\ft}$.

\subsection{Instanton component of the topological model}
When $\ell=1$, the instanton gluing-data component is given by
\begin{equation}
\label{eq:InstantonGluingDataComponent}
\bar{\Gl}_{\ft_\ell}(\delta)
=
\left(\Fr(\fg_{\ft_1})\times_X\Fr(T^*X)\times
\bar M_1^{s,\natural}(S^4,\delta)\right)/(\SO(3)\times \SO(4)).
\end{equation}
Here, $M_k^{s,\natural}(S^4,\delta)$ is the moduli space of $k$-instantons
on $S^4$, framed at the south pole $s$, with mass center at the north pole, and
scale $\leq \delta$. For $k=1$ there are homeomorphisms,
\begin{align*}
M_1^{s,\natural}(S^4,\delta)
&\cong
(0,\delta]\times\SO(3),
\\
\bar M_1^{s,\natural}(S^4,\delta)
&\cong c(\SO(3)),
\end{align*}
where $c(\SO(3))$ is the cone on $\SO(3)$.

Our gluing-model \eqref{eq:GluingModel} admits an Uhlenbeck
stratification, given below when $\ell=1$: 
$$
\bar\sM_{\ft,\fs}^{\vir}
=
\sM_{\ft,\fs}^{\vir}\sqcup \left(N_{\ft,\fs}(\eps)-M_{\fs}\right)\times X
\sqcup M_{\fs}\times X.
$$
When $\ell\geq 2$, the symmetric product has its usual stratification and
one can also construct bundles $\Gl_{\ft_\ell}(\delta,\Sigma) \to \Sigma$
with instanton moduli space fibers, for each stratum
$\Sigma\subset\Sym^\ell(X)$, following the prescription of Friedman and
Morgan \cite{FrM} and developing the idea of Kotschick and Morgan
\cite{KotschickMorgan} and Mrowka \cite{MrowkaPrivate}.  The difficult part
is to assemble these local gluing data bundles into a space of global
gluing data, $\bar{\Gl}_{\ft_\ell}(\delta)$.  One essentially has,
$$
\bar{\Gl}_{\ft_\ell}(\delta) 
=
\bigcup_{\Si\subset \Sym^\ell(X)} \Gl_{\ft_\ell}(\delta,\Si),
$$
but some modifications of the spaces $\Gl_{\ft_\ell}(\delta,\Si)$ are
needed to carry out this construction.  Unlike the analogous problem for
the anti-self-dual moduli space described in \cite{KotschickMorgan}, the
gluing maps do not define transition maps because the images of the gluing
maps for the $\SO(3)$ monopole equations intersect only at the zero-locus
of the obstruction map $\bchi$.  In \cite{FL5}, we define a deformation of
the moduli space of anti-self-dual connections on $S^4$ and deformations of
the {\em splicing \/} maps so that the intersections of the images of the
deformed splicing maps define explicitly understood transition maps.  The
space $\bar{\Gl}_{\ft_\ell}(\delta)$ is defined to be the union of the
images of these deformed splicing maps.

\subsection{Link of the Seiberg-Witten stratum}
We define the link in two steps, first considering the
Seiberg-Witten component of the link of 
$M_{\fs}\times\Sym^\ell(X)\subset\bar\sM_{\ft,\fs}^{\vir}/S^1$,
$$
\bar\bL_{\ft,\fs}^{\vir,s}
:=
\left(\rd\tN_{\ft_\ell,\fs}(\eps)
\times_{\sG_{\fs}}\bar{\Gl}_{\ft_\ell}(\delta)\right)/S^1,
$$
and, second, defining the instanton component of the link of 
$M_{\fs}\times\Sym^\ell(X)\subset\bar\sM_{\ft,\fs}^{\vir}/S^1$,
$$
\bar\bL_{\ft,\fs}^{\vir,i}
:=
\left(\tN_{\ft_\ell,\fs}(\eps)
\times_{\sG_{\fs}}\rd\bar{\Gl}_{\ft_\ell}(\delta)\right)/S^1,
$$
a complex disk bundle over
$M_{\fs}\times\rd\bar{\Gl}_{\ft_\ell}(\delta)/S^1$, when $\ell=1$, and a
more general fiber bundle (that is, not a product bundle) when $\ell>1$.
When $\ell=1$, the boundary $\rd\bar{\Gl}_{\ft_\ell}(\delta)$ is defined by
restricting to instantons on $S^4$ with scale $\delta$.
We then define the link of
$M_{\fs}\times\Sym^\ell(X)$ in the virtual moduli space,
$\bar\sM_{\ft,\fs}^{\vir}/S^1$, 
$$
\bar\bL_{\ft,\fs}^{\vir}
:=
\bar\bL_{\ft,\fs}^{\vir,s} \cup \bar\bL_{\ft,\fs}^{\vir,i}.
$$
Finally, we obtain 
$$
\bar\bL_{\ft,\fs}^{\vir}
:=
\bgamma\left(\bchi^{-1}(0)\cap \bar\bL_{\ft,\fs}^{\vir}\right)
=
(\bar\sM_{\ft}/S^1)\cap\bgamma(\bar\bL_{\ft,\fs}),
$$
the link of $M_{\fs}\times\Sym^\ell(X)$ in the true moduli space
$\bar\sM_{\ft}/S^1$.

\subsection{Fiber bundle structure in Seiberg-Witten link pairings}
The virtual moduli space method reduces the computation of the link pairing
\eqref{eq:BasicPairings} to that of the pairing on the right-hand side below,
\begin{equation}
\label{eq:VirtualPairing}
\begin{aligned}
\langle
\mu_p(z)\smile\mu_c^{\delta_c}, [\bar\bL_{\ft,\fs}]
\rangle
=
\langle
\mu_p(z)\smile\mu_c^{\delta_c}\smile e, [\bar\bL_{\ft,\fs}^{\vir}]
\rangle,
\end{aligned}
\end{equation}
where $e$ is Euler class of the total obstruction bundle over
$\bar\sM_{\ft,\fs}^{\vir}$, with section $\bchi$, and $\bar\bL_{\ft,\fs} \cong
\bchi^{-1}(0)\cap\bar\bL_{\ft,\fs}^{\vir}$.

When $\ell=1$, we show in \cite{FLLevelOne}
that the pairing \eqref{eq:VirtualPairing}
with $[\bar\bL_{\ft,\fs}^{\vir}]$ is expressible in terms of
pairings with 
$$
[M_{\fs}\times \partial\bar{\Gl}_{\ft_\ell}(\delta)/S^1].
$$
In particular, we can show that the latter pairings are in turn products of
\begin{itemize}
\item
Pairings with $[M_{\fs}]$, giving multiples of the Seiberg-Witten
invariant, $SW_X(\fs)$, and
\item
Pairings with $[\partial\bar{\Gl}_{\ft_\ell}(\delta)/S^1]$. 
\end{itemize}
The pairings with $[\partial\bar{\Gl}_{\ft_\ell}(\delta)/S^1]$ depend only
the topology of the 4-manifold, $X$, and universal data.

When $\ell>1$, the construction of the space
$$
\tN_{\ft_\ell,\fs}(\delta)\times_{\sG_{\fs}\times S^1}
\bar{\Gl}_{\ft_\ell}(\delta),
$$
shows that the virtual link $\bar\bL^{\vir}_{\ft,\fs}$ admits a
fiber bundle structure
$$
\bar\bL^{\vir}_{\ft,\fs}\to M_{\fs},
$$
with fiber given by $\CC^n\times_{S^1}\bar{\Gl}_{\ft_\ell}(\delta)$
and structure group given by $\Map(X,S^1)$.
Hence, the pairing \eqref{eq:VirtualPairing} can be written as products of
\begin{itemize}
\item
Pairings with $[M_{\fs}]$, giving multiples of the Seiberg-Witten
invariant, $SW_X(\fs)$,
\item
Pairings with the fiber $\CC^n\times_{S^1}\bar{\Gl}_{\ft_\ell}(\delta)$.
\end{itemize}
In \cite{FL5}, we show how the pairings with the fiber
$\CC^n\times_{S^1}\bar{\Gl}_{\ft_\ell}(\delta)$ can be qualitatively
understood in terms of homotopy data of the \spinu structure $\ft_\ell$,
of the \spinc structure $\fs$, and of the manifold $X$.

The passage from the rough version \eqref{eq:RoughWittenFormula} to Witten's
formula \eqref{eq:WittenFormula} to the exact formula is discussed in
\cite{FLWConjecture}.

\section{Witten's conjecture and symplectic 4-manifolds}
\subsection{Gauge theory and Lefschetz fibrations}

Away from finitely many critical points, a {\em Lefschetz fibration\/}
$\pi:X\to S$ is smooth fiber bundle over a connected base, whose fibers are
closed Riemann surfaces, $\Sigma$, of given genus. Possibly after blowing
up, all symplectic 4-manifolds admit Lefschetz fibrations
\cite{AurouxFamilySymp}, \cite{DonSympAlmostCx}, \cite{DonLefschetz}.

One can ask what is the relationship between gauge theoretic invariants and
the Lefschetz fibration structure.  As a first step, one could use the
product, $X = \Sigma\times S$, as a toy model. In this situation, the
relationship between the gauge theory moduli spaces and holomorphic maps
can be explored via the following techniques:

\begin{itemize}
\item
Adiabatic limit analysis, by the scaling metric $g=g_\Sigma\oplus g_S$ as
$g_\eps = \eps^2g_\Sigma\oplus g_S$, with $\eps\to 0$, or
\item
Restriction of stable, holomorphic bundles over $\Sigma\times S$ to
$\Sigma\times\{z\}$, as $z\in S$ varies.
\end{itemize}

This leads to the following identifications:

\begin{enumerate}
\item
SO(3)-instantons over $\Sigma\times S$ are identified with holomorphic
maps, $S\to M_\Sigma$, where $M_\Sigma$ is the space of flat SO(3)
connections over $\Sigma$. This identification goes back to Dostoglou and
Salamon \cite{DostoglouSalamonSDInst}.

\item
Seiberg-Witten U(1)-monopoles over $\Sigma\times S$ 
are identified with holomorphic maps,
$S\to M_{\Sigma,d}$, where $M_{\Sigma,d}$ is the space of vortices on the line
bundle $L\to\Sigma$ with $d=\langle c_1(L),[\Sigma]\rangle$. This
identification has been outlined by Salamon \cite{SalamonSWSympFloer}.

\item
SO(3)-monopoles over $\Sigma\times S$ are identified with holomorphic maps
from $S$ to a space of non-abelian vortices over $\Sigma$.  Non-abelian
vortices of this kind have been studied by Bradlow, Garcia-Prada, and
others.
\end{enumerate}

Of course, there are many unanswered questions:

\begin{itemize}
\item
What is the relationship between the gauge theory compactifications and
the compactifications of spaces of holomorphic maps?
\item
What is the relationship between the gauge theory invariants and
Gromov-Witten invariants of spaces of flat connections or vortices over
Riemann surfaces?
\item
Can one extend the analysis for the toy model, $\Sigma\times S$, to the
case of non-trivial surface bundles or Lefschetz fibrations? 
\item
Can symplectic 4-manifolds provide additional data we can use to 
help determine Witten's formula?
\end{itemize}

We note that the relationship between Seiberg-Witten invariants and
Lefschetz fibration data bears on a closely related program due to Simon
Donaldson and Ivan Smith whose aims include giving Lefschetz-fibration
proofs of results for symplectic 4-manifolds derived via Seiberg-Witten
theory \cite{DonaldsonSmith}.

\ifx\undefined\bysame
\newcommand{\bysame}{\leavevmode\hbox to3em{\hrulefill}\,}
\fi


\begin{thebibliography}{10}

\bibitem{AurouxFamilySymp}
D.~Auroux, {\em Asymptotically holomorphic families of symplectic
  submanifolds}, Geom. Funct. Anal. {\bf 7} (1997), 971--995.

\bibitem{DonPoly}
S.~K. Donaldson, {\em Polynomial invariants for smooth four-manifolds},
  Topology {\bf 29} (1990), 257--315.

\bibitem{DonSympAlmostCx}
S.~K. Donaldson, {\em Symplectic submanifolds and almost-complex geometry}, J.
  Differential Geom. {\bf 44} (1996), 666--705.

\bibitem{DonLefschetz}
\bysame, {\em Lefschetz fibrations in symplectic geometry}, Proceedings of the
  International Congress of Mathematicians, Vol. II (Berlin, 1998), 1998,
  pp.~309--314 (electronic).

\bibitem{DK}
S.~K. Donaldson and P.~B. Kronheimer, {\em The geometry of four-manifolds},
  Oxford Univ. Press, Oxford, 1990.

\bibitem{DonaldsonSmith}
S.~K. Donaldson and I.~Smith, {\em Lefschetz pencils and the canonical class
  for symplectic 4-manifolds}, math.SG/0012067.

\bibitem{DostoglouSalamonSDInst}
S.~Dostoglou and D.~A. Salamon, {\em Self-dual instantons and holomorphic
  curves}, Ann. of Math. (2) {\bf 139} (1994), 581--640.

\bibitem{FeehanGenericMetric}
P. M. N. Feehan, {\em Generic metrics, irreducible
rank-one $\PU(2)$ monopoles, and transversality}, Comm. Anal. Geom. {\bf 8}
(2000), 905--967, math.DG/9809001.

\bibitem{FKLM}
P.~M.~N. Feehan, P.~B. Kronheimer, T.~G. Leness, and T.~S. Mrowka, {\em {PU(2)}
  monopoles and a conjecture of {M}ari{\~n}o, {M}oore, and {P}eradze}, Math.
  Res. Lett. {\bf 6} (1999), 169--182, arXiv:math.DG/9812125.

\bibitem{FL2a}
P.~M.~N. Feehan and T.~G. Leness, 
{\em {PU(2)} monopoles and links of top-level {S}eiberg-{W}itten
  moduli spaces}, J. Reine Angew. Math. {\bf 538} (2001), 57--133,
  arXiv:math.DG/0007190.

\bibitem{FL2b}
\bysame, {\em {PU(2)} monopoles. {II}: {T}op-level {S}eiberg-{W}itten moduli
  spaces and {W}itten's conjecture in low degrees}, J. Reine Angew. Math. 
  {\bf 538} (2001), 135--212, arXiv:dg-ga/9712005. 

\bibitem{FL3}
\bysame, {\em {PU(2)} monopoles. {III}: {E}xistence of gluing and obstruction
  maps}, submitted to a print journal, arXiv:math.DG/9907107.

\bibitem{FL4}
\bysame, {\em {PU(2)} monopoles. {IV}: {S}urjectivity of gluing maps}, in
  preparation.

\bibitem{FL5},
\bysame, {\em A general SO(3)-monopole cobordism formula relating Donaldson and
Seiberg-Witten invariants}, arXiv:math.DG/0203047.

\bibitem{FLLevelOne}
\bysame, {\em {PU(2)} monopoles, level-one {S}eiberg-{W}itten moduli spaces,
  and {W}itten's conjecture in low degrees}, Topology Appl., to appear,
  arXiv:math.DG/0106238.

\bibitem{FLWConjecture},
\bysame, {\em Witten's conjecture for four-manifolds of simple type},
in preparation.

\bibitem{FSStructure}
R.~Fintushel and R.~Stern, {\em Donaldson invariants of 4-manifolds with simple
  type}, J. Differential Geom. {\bf 42} (1995), 577--633.

\bibitem{FSRationalBlowDown}
\bysame, {\em Rational blowdowns of smooth $4$-manifolds}, J. Differential
  Geom. {\bf 46} (1997), 181--235, arXiv:alg-geom/9505018.

\bibitem{FrM}
R.~Friedman and J.~W. Morgan, {\em Smooth four-manifolds and complex surfaces},
  Springer, Berlin, 1994.

\bibitem{Goettsche}
L.~G{\"o}ttsche, {\em Modular forms and {D}onaldson invariants for 4-manifolds
  with {$b^+=1$}}, J. Amer. Math. Soc. {\bf 9} (1996), 827--843,
  arXiv:alg-geom/9506018.

\bibitem{KotschickMorgan}
D.~Kotschick and J.~W. Morgan, {\em {SO(3)} invariants for four-manifolds with
  {$b^+=1$}, {II}}, J. Differential Geom. {\bf 39} (1994), 433--456.

\bibitem{KMStructure}
P.~B. Kronheimer and T.~S. Mrowka, {\em Embedded surfaces and the structure of
  {D}onaldson's polynomial invariants}, J. Differential Geom. {\bf 43} (1995),
  573--734.

\bibitem{LenessWC}
T.~G. Leness, {\em Donaldson wall-crossing formulas via topology}, Forum Math.
  {\bf 11} (1999), 417--457, dg-ga/9603016.

\bibitem{MMPdg}
M.~Mari{\~n}o, G.~Moore, and G.~Peradze, {\em Four-manifold geography and
  superconformal symmetry}, Math. Res. Lett. {\bf 6} (1999), 429--437,
  arXiv:math.DG/9812042.

\bibitem{MMPhep}
\bysame, {\em Superconformal invariance and the geography of four-manifolds},
  Comm. Math. Phys. {\bf 205} (1999), 691--735, arXiv:hep-th/9812055.

\bibitem{MrowkaPrivate}
T.~S. Mrowka, private communication.

\bibitem{PTLocal}
V.~Y. Pidstrigatch and A.~N. Tyurin, {\em Localisation of {D}onaldson
  invariants along the {S}eiberg-{W}itten classes}, arXiv:dg-ga/9507004.

\bibitem{SalamonSWBook}
D.~Salamon, {\em Spin geometry and {S}eiberg-{W}itten invariants},
  Birkh{\"a}user, Boston, to appear.

\bibitem{SalamonSWSympFloer}
D.~A. Salamon, {\em Seiberg-{W}itten invariants of mapping tori, symplectic
  fixed points, and {L}efschetz numbers}, Proceedings of 6th G\"okova
  Geometry-Topology Conference, vol.~23, 1999, pp.~117--143.

\bibitem{Witten}
E.~Witten, {\em Monopoles and four-manifolds}, Math. Res. Lett. {\bf 1} (1994),
  769--796, arXiv:hep-th/9411102.

\end{thebibliography}
\end{document}